\documentclass[preprint,12pt]{elsarticle}

\usepackage{enumerate}
\usepackage{amssymb}
\usepackage{lipsum}
\usepackage[a4paper, total={6.6in, 9.8in}]{geometry}
\usepackage{booktabs}
\usepackage{amsmath,amssymb,url}
\usepackage{enumitem} 
\usepackage{graphics} 
\usepackage{array}
\usepackage{dsfont}
\usepackage[all]{xy}
\usepackage{amsthm}
\usepackage{tikz}
\usepackage{bbding}
\usepackage{graphicx}
\usepackage{mathrsfs}

\numberwithin{equation}{section}
\newtheorem{theorem}{Theorem}[section]
\newtheorem{lemma}[theorem]{Lemma}
\newtheorem{proposition}[theorem]{Proposition}
\newtheorem{corollary}[theorem]{Corollary}

\newtheorem{problem}[theorem]{Problem}

\newtheorem{definition}[theorem]{Definition}

\newtheorem{remark}[theorem]{Remark}
\newtheorem{example}[theorem]{Example}

\newcommand{\dn}{\mathord{\downarrow}\hspace{0.05em}}
\newcommand{\up}{\mathord{\uparrow}\hspace{0.05em}}

\newcommand\blfootnote[1]{%
\begingroup
\renewcommand\thefootnote{}\footnote{#1}%
\addtocounter{footnote}{-1}%
\endgroup
}

\journal{}

\begin{document}

\begin{frontmatter}



\title{A study for well-filtered spaces}


\author{Zhengmao He}
\address{School of Sciences, Southwest Petroleum University,Chengdu 610500,Sichuan China}

\begin{abstract} In this paper, we focus on well-filtered spaces. By using the Rudin sets and WD sets, we obtain the following results:

(1) The product of arbitrarily many Scott well-filtered dcpos is still well-filtered with respect to the Scott topology. Thus we give positive answers for two questions proposed by Xu and Zhao;

(2) Let $P$ be a dcpo with $\Sigma K(P)$ being well-filtered and $Q$ an RB-domain. Then $\Sigma K(P\times Q)$ is well-filtered if $K(P\times Q)$ is a dcpo;

(3) An intrinsic characterization of the $WD$ sets is given in terms of Scott open filters.
\end{abstract}
\begin{keyword}  Scott topology; well-filtered; RB domain; WD set\\
\vspace*{0.2cm}
\end{keyword}


\end{frontmatter}
\blfootnote{This work is supported by the National Natural Science Foundation of China (Grant nos.12471438,12601906) and the Sichuan Science and Technology Program (Grant No. 2026NSFSC0784).}
\blfootnote{E-mail address: hezhengmaomath@163.com (Z.M.He).}


\section{Introduction}
\label{}

Domain theory constitutes a framework within which order and topology admit a profound integration. The Scott topology is a canonical topology induced by the order structure of a poset and plays a fundamental role in domain theory. Sobriety and well-filteredness are two kinds of important topological properties. The study of sobriety and well-filteredness for the Scott spaces clarifies two complementary aspects of the interaction between order and topology. The Scott topology on a continuous dcpo possesses favourable properties. If we extend
these properties to general dcpos, it requires some topological conditions that
remain useful without Hausdorff separation axiom. Well-filteredness is a desirable candidate. In this paper, we focus on well-filtered spaces.

The distinction between well-filteredness and sobriety emerged from
a series of counterexamples and representation results. Johnstone
constructed a dcpo whose Scott topology is not sober and not well-filtered(\cite{Johnstone81}). As we know,
every sober space is well-filtered(\cite{GG03,JG13}). Kou showed
that there is a dcpo which is well-filtered but non sober for the Scott topology and thus answered a question posed by Heckmann (\cite{Kou01,Heckmann92}.
Xi and Zhao later study the models for well-filtered $T_1$ spaces
and produce new non sober examples (\cite{XiZhao17}).
Xi and Lawson proved that every Lawson-compact dcpo
and hence every complete lattice is Scott well-filtered(\cite{XL17}). In addition, Lawson, Wu and Xi proved that every
core-compact well-filtered $T_0$ space is sober (\cite{LWX20}).
Topological Rudin's lemma extended by Heckmann
and Keimel provides a useful tool for constructing an irreducible closed set(\cite{HK13}). Shen, Xi, Xu and Zhao constructed well-filtered reflections and proved that arbitrary ordinary topological products of well-filtered spaces are
well-filtered(\cite{SXXZ19}). Furthermore, Xu, Shen, Xi and Zhao
developed Rudin
sets and well-filtered determined sets( or $\mathrm{WD}$ sets) to
characterize well-filteredness and describe its reflections(\cite{XSXZ20}). Thus the Rudin sets and $\mathrm{WD}$ sets are useful tools for studying well-filtered spaces.

The product preservation of well-filteredness and sobriety is already nontrivial for the Scott topology. For two posets $P$ and $Q$, the Scott topology
$\Sigma(P\times Q)$ can be strictly finer than the product topology
$\Sigma P\times\Sigma Q$. Consequently, the ordinary product theorems with respect to sobriety and well-filteredness do not apply directly to this Scott product problem. So Xu and Zhao
propose four questions for well-filtered spaces and sober spaces (Questions 4.6,4.7,4.8 and 4.9 in\cite{XZ21}).
Recently, Miao, Xi,Jia, Li and Zhao construct two Scott sober
dcpos whose order product is not Scott sober(\cite{MXJLZ26}). Later, He obtained two Scott sober complete lattices even two co-frames whose order product is not sober for the Scott topology(\cite{HMOOP281}). In this paper, we solve the question 4.7 and  question 4.9. Specifically, the first result is
that for every family of Scott well-filtered dcpos $\{P_i\mid i\in I\}$,
$\Sigma\prod\limits_{i\in I}P_i$ is well-filtered.
The proof for the above result analyzes projections of
Rudin sets and finite coordinate replacements.

The collection $K(X)$ consists of all the nonempty compact saturated subsets of $X$ ordered by the reverse inclusion. The upper Vietoris topology
of the Smyth power space and the Scott topology of $K(X)$ need not
coincide. Xu, Wen and Xi proved that well-filteredness
of $X$ implies that  $\Sigma K(X)$ is well-filtered(\cite{XWX23}). However, for two dcpos $P$ and $Q$, $K(P\times Q)$ is not equal to $K(P)\times K(Q)$. So the well-filteredness of $\Sigma K(P\times Q)$ is unknown under the condition that $\Sigma K(P)$ and $\Sigma K(Q)$ are well-filtered. We first prove that for dcpos $P$ and $Q$, $\Sigma K(P)$ and $\Sigma K(Q)$ are well-filtered if
$\Sigma K(P\times Q)$ is well-filtered. Furthermore, we also show that if $\Sigma K(P)$ is well-filtered and $Q$ a $RB$-domain such that $K(P\times Q)$ is a dcpo, then $\Sigma K(P\times Q)$ is well-filtered.

The $\mathrm{WD}$ sets as a useful tool for studying well-filtered sapces are defined by using all well-filtered spaces and continuous mappings. Thus it is difficult to check whether a closed subset of a topological space is a $\mathrm{WD}$ set. It is essential to give a characterization for $\mathrm{WD}$ sets of a topological space $X$ by just using the information on $X$.  Motivated by the open filters viewpoint of the
Hofmann--Mislove theorem, we prove that, for every closed
$A\subseteq X$,
$A\in\mathrm{WD}(X)$ if and only if $p_A\in\mathscr H_X$, where $\mathscr H_X$ is a subfamily of all proper Scott open filters on $\mathcal{O}(X)$ with some rules.
\section{Preliminaries}
\quad Let $P$ be a poset. For each $A\subseteq P$, set $\up A=\{x\in P\mid \exists\ a\in A, a\leq x \}$ and $\dn A=\{x\in P\mid \exists\ a\in A, x\leq a\}$. We write $\dn\{x\}$ as $\dn x$ and $\up\{x\}$ as $\up x$ for each $x\in P$. A subset $D$ of $P$ is {\em directed} provided that $D$ is nonempty and $\forall\ a,b\in D$, $\up a\cap \up b \cap D\neq\emptyset$. Dually, a nonempty subset $D\subseteq P$ is {\em filtered} if $\forall\ a,b\in D$, $\dn a\cap \dn b \cap D\neq\emptyset$. An element $x$ is an {\em upper bound} for a subset $A$ if $a\leq x$ for each $a\in A$. A poset $P$ is called a {\em directed complete poset} ({\em dcpo}, for short) if every directed subset of $P$ has a supremum. For each pair of $x, y\in P$, we say that $x$ is {\em way below} $y$, written $x\ll y$, if for each directed set $D\subseteq P$ whose supremum exists, $y\leq\bigvee D$ implies $\up x\cap D\neq\emptyset$. An element $p\in P$ is {\em compact} if $p\ll p$. A subset $U\subseteq P$ is {\em Scott open} if $U=\up U$ and for each directed subset $D$ with $\bigvee D$ existing, $\bigvee D\in U$ implies $D\cap U\neq \emptyset$. The set of Scott open sets of $P$ forms the Scott topology $\sigma(P)$. We write $\Sigma P$ for $(P,\sigma(P))$. Let $P,Q$ be dcpos. Then the mapping $f:\Sigma P\longrightarrow \Sigma Q$ is continuous (or called Scott continuous) if and only if $f$ is monotone and for each directed set $D\subseteq P$, $f(\bigvee D)=\bigvee\limits_{d\in D}f(d)$.

\quad Let $X$ be a topological space. We write \(\operatorname{cl}_{X}(A)\) for the closure of $A$. A subset $A$ of $X$ is called a {\em saturated set} if $A$ is just the intersection of all open sets containing $A$. A subset $F$ of $X$ is {\em irreducible} if $F$ is nonempty and for every pair of closed sets $A, B\subseteq X$, $F\subseteq A\cup B$ implies $F\subseteq A$ or $F\subseteq B$. A $T_{0}$ space $X$ is {\em sober} if every irreducible closed set is just the closure of a singleton. For a $T_{0}$ space $X$, the {\em specialization order} $\leq$ on $X$ is defined by $x\leq y $ if and only if $x\in cl_{X}(\{y\})$. Unless otherwise stated, throughout
the paper, whenever an order-theoretic concept is mentioned in the context of a $T_{0}$ space $X$, it is to be interpreted with respect to the specialization order on $X$. A subset $A$ of a topological space $X$ is {\em saturated} if $A=\up A$ under the specialization order. We use $K(X)$ to denote the set of all nonempty compact saturated sets of $X$ with the reverse inclusion order and $\mathcal{O}(X)$ ($\Gamma(X)$) to denote the lattice of all open subsets (closed subsets) of $X$ ordered by inclusion. A topological space $X$ is {\em well-filtered} if for each filtered family $\mathcal{F}$ of compact saturated subsets of $X$ and each open set $U$ of $X$, $\bigcap \mathcal{F}\subseteq U$ implies $F\subseteq U$ for some $F\in\mathcal{F}$. A dcpo $P$ is called {\em Scott well-filtered} if the dcpo $P$ endowed with the Scott topology is well-filtered. A topological space $X$ is a {\em retract} of a topological space $Y$ if there are two continuous maps $f:X\longrightarrow Y$ and $g:Y\longrightarrow X$ such that $g\circ f=id_{X}$.

\section{The Scott well-filteredness of the arbitrary product}

In this section, we show that the product of Scott well-filtered dcpos is still Scott well-filtered. Thus this gives positive answers to Questions 4.7 and 4.9 posed by Xu and Zhao(\cite{XZ21}).

\begin{definition}
{\rm (\cite{XSXZ20}) Let $X$ be a $T_0$ space. A nonempty closed subset $A$ of $X$ is called a
\emph{Rudin set} if there exists a filtered family
$\mathscr{K}\subseteq K(X)$ such that $A$ is minimal in
$$M(\mathscr{K})=\{B\in\Gamma(X)\mid \forall\ K\in\mathscr{K}, B\cap K\neq\emptyset\},$$
where the order on $M(\mathscr{K})$ is the inclusion order. The symbol ${\rm RD}(X)$ denotes the collection of all closed Rudin sets of $X$.}
\end{definition}

In the following, we call the filtered family $\mathscr{K}\subseteq K(X)$ a witness of the set $A\in{\rm RD}(X)$.

\begin{lemma}{\rm
Let $\mathscr{K}$ be a witnesses of a closed Rudin set $A$. Then for every $K_0\in\mathscr{K}$,
$$
  cl_X(A\cap K_0)=A.
 $$
}
\end{lemma}

\begin{proof}
Choose any $K\in\mathscr{K}$. By the filteredness of $\mathscr{K}$, there exists  $K_{1}\in\mathscr{K}$ with $K_{1}\subseteq K_0\cap K$. It follows from $A\in M(\mathscr{K})$ that we can select $a\in A\cap K_{1}$. This implies $a\in A\cap K_0\cap K$ and thus
$$
  K\cap cl_X(A\cap K_0)\ne\varnothing.
$$
So $cl_X(A\cap K_0)\in M(\mathscr{K})$ and clearly $cl_X(A\cap K_0)\subseteq A$. The minimality of $A$ gives $cl_X(A\cap K_0)=A$.
\end{proof}

\begin{theorem}{\rm (\cite{XSXZ20})
Let $X$ be a $T_0$ space. Then the following are equivalent:

(1) $X$ is well-filtered;

(2) For every closed Rudin set $A$, $A=cl_X(\{a\})$ for some $a\in X$.}
\end{theorem}

\begin{lemma}
{\rm Let $R=\prod\limits_{i\in I}P_i$ be the product of a family consisting of dcpos and
$A\in {\rm RD}(\Sigma R)$. For each $i\in I$, set $A_i=\pi_i(A)$.
Then the following statements hold.

(1) $A_i\in {\rm RD}(\Sigma P_i)$;

(2) If $a=(a_{i})_{i\in I}\in A$ and $b_j\in A_j$, then $x=[a,b_{j}]\in A$, where $x=(x_{i})_{i\in I}$ with $x_{i}=a_{i}$ ($i\in I\setminus\{j\}$) and $x_{j}=b_{j}$ ;

(3) If $a\in A$ and for any finite $F\subseteq I$, $b_j\in A_j$ for $j\in F$, then $y=[a,F]\in A$, where $y=(y_{i})_{i\in I}$ with $y_{i}=a_{i}$ ($i\in I\setminus F$) and $y_{j}=b_{j}$ ($j\in F$).}

\end{lemma}

\begin{proof}
(1) As $A\in{\rm RD}(\Sigma P)$, there exists a filtered family $\mathscr{K}=\{K_{\lambda}\mid \lambda\in\Lambda\}$ consisting of nonempty compact saturated subsets of $\Sigma R$ such that $A$ is a minimal element in $M(\mathscr{K})$.

Fix $i\in I$ and set
$$
 Q=\prod_{j\in I\setminus\{i\}}P_j,\qquad T=P_i.
$$
Then $Q$ and $T$ are dcpos. Define a mapping $\varphi_i:R\longrightarrow Q\times T$ by
$$\varphi_i(x)=(x|_{I\setminus\{i\}},x_i),$$ where $x|_{I\setminus{i}}$ denotes the element obtained from $x$ by deleting its $i$-th coordinate. One can check that $\varphi_i$ is an order isomorphism. Moreover, $\varphi_i$ preserves the directed suprema coordinatewise. Consequently, $\varphi_i$ and its converse mapping are Scott continuous.
Set $B=\varphi_i(A)$ and let
 $$\pi_{Q}:Q\times T\longrightarrow Q, \pi_{T}:Q\times T\longrightarrow T
$$
be two projections. Thus $\pi_{T}(B)=A_i$. Furthermore, $\pi_{Q}$ and $\pi_{T}$ are Scott continuous. Take
$$
 L_\lambda=\varphi_i(K_\lambda).
$$
Each $ L_\lambda$ is a nonempty saturated compact set in $\Sigma(Q\times T)$ because $\varphi_i$ is an order isomorphism. The collection $\{L_{\lambda}\mid \lambda\in\Lambda\}$ is a filtered family as $\{\lambda\mid \lambda\in\Lambda\}$ is filtered.

$\mathbf{Claim \ 1}$: $\{L_{\lambda}\mid \lambda\in\Lambda\}$ is a witness family for the closed set $B$.

Since $\varphi_i$ is an order isomorphism and $A$ is Scott closed in $R$, $B$ is closed in $\Sigma (Q\times T)$.
Suppose $C\subseteq B$ is a Scott closed subset of $Q\times T$ which meets each $L_\lambda$. Then
$\varphi_i^{-1}(C)\subseteq A$ is a closed subset of $R$ and meets each $K_\lambda$.
The minimality of $A$ implies that $\varphi_i^{-1}(C)=A$. As a consequence, $C=B$ and hence $B\in{\rm RD}(\Sigma (Q\times T))$. This proves Claim 1.

For each $\lambda$, set
$$ S_\lambda=B\cap L_\lambda,\ H_\lambda=\pi_{Q}(S_\lambda)\text{ and } J_\lambda=\pi_{T}(S_\lambda).$$

As $B$ is closed in $\Sigma(Q\times T)$ and $L_{\lambda}$ is compact, $S_{\lambda}$ is also compact. By the Scott continuity of $\pi_{Q}$ and $\pi_{T}$, $H_\lambda$ and $J_\lambda$ are nonempty compact subset of $\Sigma Q$ and $\Sigma T$, respectively. By Lemma 3.2, we have
$$
  cl_{\Sigma(Q\times T)} (S_\lambda)=B.
$$
For each $\lambda$, define
$$
 E_\lambda
 =\{t\in T\mid \exists\ q\in H_\lambda,\ (q,t)\in B\}.
$$

$\mathbf{Claim\ 2}$: $E_\lambda$ is closed in $\Sigma T$.

Let $t'\leq t\in E_\lambda$ and $(q,t)\in B$ with $q\in H_\lambda$. Then
$(q,t')\leq (q,t)$ and further $(q,t')\in B$ because $B$ is lower. So $t'\in E_\lambda$ and therefore $E_\lambda$ is also lower.

Suppose $D\subseteq E_\lambda$ is directed and take
$\hat{t}=\sup\limits_{T} D$. For each $d\in D$, put
$$
 F_d=\{q\in H_\lambda:(q,d)\in B\}.
$$
Then $F_d$ is nonempty as $d\in E_\lambda$. It is not difficult to check that the mapping $f_{d}:Q\longrightarrow Q\times T$ is Scott continuous, where $f_{d}$ is defined by $$\forall a\in Q, f_{d}(a)=(a,d).$$ Consequently, $$f_{d}^{-1}(B)=\{q\in Q\mid (q,d)\in B\}$$ is Scott closed in $Q$.
So $F_d=H_{\lambda}\cap f_{d}^{-1}(B)$ is closed in the compact subspace $H_\lambda$.
Choose any $d_1,\ldots,d_n\in D$. There is a $d_0\in D$ larger than each $d_{k}$ ($1\leq k\leq n$) as $D$ is directed. Since $B$ is lower, we have
$$
 \varnothing\ne F_{d_0}\subseteq\bigcap_{k=1}^n F_{d_k}.
$$
The compactness of $H_{\lambda}$ gives that we can select $\hat{q}\in\bigcap\limits_{d\in D}F_d\subseteq H_{\lambda}$. Now,
$\{(\hat{q},d)\mid d\in D\}\subseteq B$ is directed and $$\bigvee\limits_{d\in D}(\hat{q},d)=(\hat{q},\hat{t}).$$ It follows from the Scott closedness of $B$ that $(\hat{q},\hat{t})\in B$. Then we deduce that $\hat{t}\in E_\lambda$. Hence, $E_\lambda$ is Scott closed in $T$.

$\mathbf{Claim\ 3}$: Each $A_{i}$ is closed in $\Sigma T$.

By the definition of these sets, we can check that
$$
 J_\lambda\subseteq E_\lambda\subseteq A_i.
$$
Then $ cl_{\Sigma T}(J_\lambda)\subseteq cl_{\Sigma T}(A_i)$.

As $$S_{\lambda}\subseteq\pi_{T}^{-1}(\pi_{T}(S_{\lambda}))\subseteq\pi_{T}^{-1}( cl_{\Sigma T}(J_\lambda)),$$ we know $$B=cl_{\Sigma(Q\times T)}(S_{\lambda})\subseteq\pi_{T}^{-1}( cl_{\Sigma T}(J_\lambda)).$$
Hence $A_{i}=\pi_{T}(B)\subseteq cl_{\Sigma T}(J_\lambda)$. This shows that $ cl_{\Sigma T}(J_\lambda)=cl_{\Sigma T}(A_i)$.

By Claim 2, $E_\lambda$ is Scott closed in $T$. Consequently, we have
$$
  cl_{\Sigma T}(A_i)
 = cl_{\Sigma T}(J_\lambda)
 \subseteq E_\lambda
 \subseteq A_i
 \subseteq cl_{\Sigma T}(A_i).
$$
This yields that for every $\lambda$,
$$
 E_\lambda=A_i= cl_{\Sigma T}(A_i).
$$
This proves Claim 3.

$\mathbf{Claim\ 4}$: Each $A_i$ is a closed Rudin set in $\Sigma T$.

Set
$$
 M_\lambda=\up_T J_\lambda
 =\{t\in T\mid \exists\ s\in S_\lambda,\pi_{T}(s)\le t\}.
$$
Then each $M_\lambda$ is a nonempty compact saturated set as $J_{\lambda}$ is compact.

If $L_\nu\subseteq L_\lambda\cap L_\mu$, then
$$
 S_\nu\subseteq S_\lambda\cap S_\mu\text{ and thus }M_\nu\subseteq M_\lambda\cap M_\mu.
$$
So the family $\mathscr{A}=\{M_\lambda\mid \lambda\in\Lambda\}$ is a filtered family of nonempty compact saturated set in $\Sigma T$. Moreover, for each $\lambda\in\Lambda$,
$$A_i\cap M_\lambda\supseteq J_\lambda\neq\emptyset.$$
In other words, the closed set $A_{i}$ belongs to $M(\mathscr{A})$.

Let $C\subseteq A_i$ be a closed subset in $\Sigma T$ such that
$C$ meets each $M_\lambda$. Choose $c\in C\cap M_\lambda$. Then there is
$s\in S_\lambda$ with $\pi_{T}(s)\leq c$. As $C$ is lower,
we obtain $\pi_{T}(s)\in C$. As a result,
$$
 s\in B\cap\pi_{T}^{-1}(C)\cap L_\lambda.
$$
This implies that $B\cap\pi_{T}^{-1}(C)$ is a closed subset of $B$ that meets each
$L_\lambda$. By Claim 1, the minimality of $B$ gives $B\cap\pi_{T}^{-1}(C)=B$.
This yields that $A_i=\pi_{T}(B)\subseteq C$ and therefore $C=A_i$. This proves Claim 4.

(2) Define a mapping $\pi_{-i}:R\longrightarrow Q$ by $\pi_{-i}(r)=r|_{I\setminus\{i\}}$.

$\mathbf{Claim\ 5}$:  $\varphi_i(A)=\pi_{Q}(B)\times A_i=\pi_{-i}(A)\times A_i$.

Fix $t\in A_i$ and set
$$
 B_t=\{q\in Q\mid(q,t)\in B\}.
$$
Then $B_{t}=f_{t}^{-1}(B)$ is Scott closed in $Q$ as $f_{t}$ is Scott continuous. As for every $\lambda$, $t\in A_i=E_\lambda$, there
exists $q_\lambda\in H_\lambda$ satisfying $(q_\lambda,t)\in B$.
So $q_\lambda\in B_t$ and besides there
is $u_\lambda\in T$ with
$$
 (q_\lambda,u_\lambda)\in S_\lambda=B\cap L_\lambda.
$$
So we know for each $\lambda\in\Lambda$,
$$
(q_\lambda,u_\lambda)\in B\cap\pi_{Q}^{-1}(B_t)\cap L_{\lambda}.$$

Now, the closed subset $B\cap\pi_{Q}^{-1}(B_t)\subseteq B$ meets each $L_\lambda$. By Claim 1, $B\cap\pi_{Q}^{-1}(B_t)=B$ and hence each $q\in\pi_{Q}(B)$ belongs to $B_t$.
As a consequence, $(q,t)\in B$ for every $q\in\pi_{Q}(B)$.

Since $t\in A_i$ is arbitrary, this indicates that $\pi_{Q}(B)\times A_{i}\subseteq B$. It is not difficult to verify that $B\subseteq\pi_{Q}(B)\times A_{i}$ and $\pi_{Q}(B)=\pi_{-i}(A)$. Therefore, we have
$$
 \varphi_i(A)=\pi_{Q}(B)\times A_i
             =\pi_{-i}(A)\times A_i.
$$
This Claim is proved.

$\mathbf{Claim\ 6}$: Let $a\in A$ and $b_i\in A_i$. Then $[a,b_{i}]\in A$, where $[a,b_{i}]\in R$ is defined by
$$
 \bigl([a,b_{i}]\bigr)_j=
 \begin{cases}
 b_i,&j=i,\\
 a_j,&j\neq i.
 \end{cases}
$$
As $\pi_{-i}(a)\in\pi_{-i}(A)$, by Cliam 5, we get
$$
 \varphi_i\bigl([a,b_{i}])
       =(\pi_{-i}(a),b_i)\in\varphi_i(A).
$$
Then $[a,b_{i}]\in A$ because $\varphi_i$ is an order isomorphism. The assertion (2) is valid.

(3) Fix $a\in A$ and a finite set $F\subseteq I$ with $b_i\in A_i$ for each
$i\in F$. Define
$$
 [a,F]_j=
 \begin{cases}
 b_j,&j\in F,\\
 a_j,&j\notin F.
 \end{cases}
$$
$F=\varnothing$ means $[a,F]=a$. So, we consider the situation $F\neq\emptyset$ and
enumerate $F$ as $\{i_1,\ldots,i_n\}$. Set
$$
 a^{(0)}=a,\
 a^{(k)}=[ a^{(k-1)},b_{i_{k}}]
 \quad(1\le k\le n).
$$
By induction, each $a^{(k)}$ ($0\leq k\leq n$) belongs to $A$ and hence $[a,F]=a^{(n)}\in A$. The assertion (3) holds.
\end{proof}

\begin{theorem}
{\rm Let $\{P_i\mid i\in I\}$ be a family of dcpos such that ecah $\Sigma P_i$
is well-filtered. Then $\Sigma\prod\limits_{i\in I}P_i$ is well-filtered.}
\end{theorem}

\begin{proof}
By Theorem 3.3, it suffices to prove that each
$A\in {\rm RD}(\Sigma R)$ is a principal ideal.

By Lemma 3.4, $A_i=\pi_i(A)$ is a closed
Rudin set in $\Sigma P_i$. By Theorem 3.3, there is a unique
$p_i\in P_i$ such that $A_i=\dn p_i.$
Set $p=(p_i)_{i\in I}\in R$. Then for each $a\in A$ and $i\in I$,
$a_i\in A_i=\dn p_i$. As a result, $a\in \dn p$. This yields that
$p$ is larger than any elements coming from $A$. It remains to show that $p\in A$.

As $A\ne\varnothing$, we can select a $b=(b_i)_{i\in I}\in A$.
Let ${\rm Fin}(I)$ be the collection of all finite subsets of $I$, containing the empty set. For each $F\in {\rm Fin}(I)$, define
$$
 b_i^F=
 \begin{cases}
   p_i,&i\in F,\\
   b_i,&i\notin F.
 \end{cases}
$$
As each $p_i\in A_i$ and $b\in A$, by Lemma 3.4(3),
we know $b^{F}=(b_{i}^{F})_{i\in I}\in A$ for every $F\in{\rm Fin}(I)$. So each $b^F\leq p$.
Since each $p_{i}$ is the largest element in $A_{i}$, we deduce that $b^F\leq b^G$ if $F\subseteq G$. Thus the set
$$
D=\{b^F\mid F\in Fin(I)\}
$$
is directed and contained in $A$. Clearly, $p$ is an upper bound of $D$. For each $i\in I$, choose $\widetilde{F}=\{i\}$. Then $b^{\widetilde{F}}\in D$ and hence $p_{i}=\pi_{i}(b^{\widetilde{F}})\in\pi_{i}(D)$. Consequently, we obtain
$$
  \sup\limits_{R} D=(\bigvee\limits_{P_{i}}\pi_{i}(D))_{i\in I}=(p_{i})_{i\in I}=p.
$$

Since $A$ is Scott closed and $D\subseteq A$ is directed, we have
$p=\sup D\in A$. Now, we can conclude that $p$ is a largest element of $A$. Therefore,
$A=\dn p$. By Theorem 3.3, $\Sigma R$ is well-filtered.
\end{proof}

\begin{lemma}{\rm (\cite{XiZhao17}) A retract of a well-filtered space is well-filtered.}

\end{lemma}
\begin{corollary}
{\rm Let $R=\prod_{i\in I}P_i$ be a nonempty order product of dcpos. Then the following statements are equivalent:

(1) $\Sigma R$ is well-filtered;

(2) Each $\Sigma P_i$ is well-filtered.}
\end{corollary}

\begin{proof}

$(1)\Longrightarrow(2)$ Since each $\Sigma P_{i}$ is a retract of $\Sigma R$, by Lemma 3.6, we have $\Sigma P_i$ is well-filtered.

$(2)\Longrightarrow(1)$ It is straightforward by using Theorem 3.5.
\end{proof}

\section{The Scott well-filteredness of the product for Smyth powerdomain}

In this section, we consider when $\Sigma K(P\times Q)$ is well-filtered under the assumption that $\Sigma K(P)$ and $\Sigma K(Q)$ are well-filtered.

\begin{lemma}{\rm
Let $P,Q$ be a pair of dcpos and $H\in K(P)$. If
$C\subseteq P\times Q$ is Scott closed, then $C(H)$ is Scott closed in $Q$, where
$$ C(H)=\pi_Q(C\cap(H\times Q))=\{q\in Q\mid \exists p\in H,\ (p,q)\in C\}
.$$
}
\end{lemma}

\begin{proof}

$\mathbf{Claim\ 1}$: $C(H)$ is a lower set.

Suppose $q\in C(H)$ and $q_{0}\leq q$. Then there exists $p\in H$ with $(p,q)\in C$.
As $C\subseteq P\times Q$ is a lower set and $(p,q_{0})\leq(p,q)$, we have
$(p,q_{0})\in C$. As a result, $q_{0}\in C(H)$. This proves Claim 1.

$\mathbf{Claim\ 2}$: $C(H)$ is closed under directed suprema in $Q$.

Let $D\subseteq C(H)$ be a directed subset and  $q_{1}=\sup D$.
For every $d\in D$, take
$$
  F_d=\{h\in H\mid (h,d)\in C\}=H\cap C_d,
$$
where $C_{d}=\{x\in P\mid (x,d)\in C\}$.
Then $F_d$ is nonempty as $d\in C(H)$. One can check that $C_d$ is
Scott closed in $P$. Choose any $d_1,\ldots,d_n\in D$. The directedness of $D$ provides $d_0\in D$ with $d_k\leq d_0$ for each $k$. If $p\in F_{d_0}$, then $(p,d_0)\in C$.
Since $C$ is a lower set, each $(p,d_k)$ belongs to $C$. In other words,
$$
  \varnothing\ne F_{d_0}\subseteq\bigcap_{k=1}^n F_{d_k}.
$$
As $H$ is compact, we can choose a
$p_{0}\in\bigcap_{d\in D}F_d$. Equivalently, $(p_0,d)\in C$ for each $d\in D$.

Now, the directed set $\{(p_0,d)\mid d\in D\}$ is contained in $C$ and further $$\bigvee\limits_{d\in D}(p_0,d)=(p_0,q_{1}).$$
It follows from the Scott closedness of $C$ that $(p_0,q_{1})\in C$.
As $p_0\in H$, this yields that $q_1\in C[H]$. This proves Claim 2.

By these two claims, $C(H)$ is Scott closed in $Q$.
\end{proof}

\begin{remark}
{\rm In Lemma 4.1,it suffices to assume that $H$ is a compact subset of $\Sigma P$. The saturation of $H$ is useless.}
\end{remark}

\begin{corollary}
{\rm Let $H$ be a compact subset in $\Sigma P$ and $U\in\sigma(P\times Q)$. Then
$$
 V_H(U)=\{q\in Q\mid H\times\{q\}\subseteq U\}
$$
is Scott open in $Q$.}
\end{corollary}

\begin{proof}
Set $C=(P\times Q)\setminus U$. One can check that $$
 Q\setminus V_H(U)
   =\{q\in Q\mid \exists p\in H,\ (p,q)\in C\}=C(H).
$$
By Lemma 4.1 and Remark 4.2 , $C(H)$ is Scott closed in $Q$. Hence, $V_H(U)\in\sigma(Q)$.
\end{proof}

\begin{proposition}
{\rm Let $H$ be a compact subset in $\Sigma P$ and $K$ be a compact subset in $\Sigma Q$.
Then $H\times K$ is compact in $\Sigma(P\times Q)$. Furthermore,
if $H,K$ are nonempty and saturated, then $H\times K$ is also nonempty
compact and saturated in $\Sigma(P\times Q)$.}
\end{proposition}

\begin{proof}
This proposition holds if either $H$ or $K$ is empty. So we suppose that
$H$ and $K$ are nonempty. Let $\{U_i\mid i\in I\}\subseteq\sigma(P\times Q)$ be an open
cover of $H\times K$. For every finite subset $F\subseteq I$, set
$$
 U_F=\bigcup_{i\in F}U_i\text{ and }
 V_F=\{q\in Q\mid H\times\{q\}\subseteq U_F\}.
$$
By Corollary 4.3, $V_F\in\sigma(Q)$. Clearly,
$V_F\subseteq V_G$ if $F\subseteq G$.

Fix $q\in K$. For each $i\in I$, it follows from $U_{i}\in\sigma(P\times Q)$ that the set
$$
 W_{i,q}=\{p\in P\mid (p,q)\in U_i\}
$$
is Scott open in $\Sigma P$. As all $U_{i}$ cover $H\times K$, the collection $\{W_{i,q}\mid i\in I\}$ covers the compact set $H$ for each $q\in K$. The compactness of $H$ gives that there is a finite $F_{q}\subseteq I$
such that $$H\subseteq\bigcup\limits_{i\in F_{q}}W_{i,q}.$$ Equivalently, $H\times\{q\}\subseteq U_{F_{q}}$ and therefore $q\in V_{F_{q}}$.
Consequently, $\{V_F\mid F\subseteq I\text{ finite}\}$ forms an open cover
of $K$. By the compactness of $K$, there exist finite subsets
$F_1,\ldots,F_m$ contained in $I$ satisfying $$K\subseteq V_{F_1}\cup\cdots\cup V_{F_m}.$$
Put $F^{\ast}=F_1\cup\cdots\cup F_m$. Then $F^{\ast}$ is finite and satisfies
$$
  K\subseteq V_{F^{\ast}}\text{ and } H\times K\subseteq U_{F^{\ast}}.
$$
This shows that $H\times K$ is compact in $\Sigma(P\times Q)$. The saturation of $H\times K$ is valid if $H$ and $K$ are saturated.
\end{proof}

By Proposition 4.4, for two dcpos $P$ and $Q$, $K(P)\times K(Q)$ can be order embedded into $K(P\times Q)$ by the product operator. The following example shows that $K(P)\times K(Q)$ is not equal to $K(P\times Q)$ in general.

\begin{example} {\rm Let $P=\mathbb{N}\cup\{\infty\}$ be the dcpo ordered by $$0\leq 1\leq\cdot\cdot\cdot\leq n\leq n+1\leq\cdot\cdot\cdot\leq\infty.$$
Then $P$ is a complete lattice and hence $\Sigma P$ is well-filtered. One can verify that $K(P)=\{\up x\mid x\in P\}$. Consider $A=\up_{P\times P}(1,0)\cup\up_{P\times P}(0,1)$. Clearly, $A\in K(P\times P)$ but $A\not\in K(P)\times K(P)$. So we have $K(P)\times K(P)\neq K(P\times P)$. Furthermore, $K(P\times P)$ is not isomorphic to $K(P)\times K(P)$. Note that $K(P)$ is order isomorphic to $P$. Choose any three elements $\alpha=(a,b)$, $\beta=(c,d)$ and $\gamma=(e,f)$ in $P\times P$. Then we have
$$
\begin{aligned}
&\quad \alpha\wedge \beta\wedge \gamma\\
&\quad=({\rm min}\{a,c,e\},{\rm min}\{b,d,f\})\\
&\quad={\rm min}\{\alpha\wedge \beta,\alpha\wedge\gamma,\beta\wedge \gamma\}.
\end{aligned}
$$
Set $B=\up_{P\times P}(0,2)$, $C=\up_{P\times P}(1,1)$ and $D=\up_{P\times P}(2,0)$ in $K(P\times P)$. Then $B,C$ and $D$ belong to $K(P\times P)$ and $$B\wedge C\wedge D=B\cup C\cup D\not\in\{B\wedge C,B\wedge D,C\wedge D\}.$$ This means that $K(P\times P)\not\cong P\times P$. Therefore, $K(P\times P)\not\cong K(P)\times K(P)$.

}
\end{example}

\begin{lemma} {\rm (\cite{GG03}) Let $X$ be a $T_{0}$ space and $\mathscr{A}\subseteq K(X)$ a nonempty family. Then $\bigvee\limits_{K(X)}\mathscr{A}$ exists if and only if $\bigcap\mathscr{A}\in K(X)$. In this case, $\bigvee\limits_{K(X)}\mathscr{A}=\bigcap\mathscr{A}$.}
\end{lemma}

\begin{lemma}
{\rm (\cite{GG03}) Let $X$ be a well-filtered $T_{0}$ space and $\mathcal{K}$ a filtered family consisting of nonempty compact saturated subsets of $X$. Then $\bigcap\mathcal{K}$ is also a nonempty compact  saturated set.}
\end{lemma}

\begin{theorem} {\rm Let $P$ and $Q$ be dcpos. If $\Sigma K(P\times Q)$ is well-filtered, then $\Sigma K(P)$ and $\Sigma K(Q)$ are well-filtered.}

\end{theorem}

\begin{proof} We first prove that $K(P)$ and $K(Q)$ are dcpos.
Fix a $q_0\in Q$, define a mapping $e_P:K(P)\longrightarrow K(P\times Q)$ by $e_P(A)=A\times\uparrow_Qq_0$. By Proposition 4.4, $e_{P}$ is well defined.

Let $\{A_i\mid i\in I\}$ be a nonempty
filtered family in $K(P)$. Then
$\{e_P(A_i)\mid i\in I\}$ forms a directed subset of
$K(P\times Q)$. Since $\Sigma K(P\times Q)$ is well-filtered, $K(P\times Q)$ is a dcpo. By Lemma 4.6 and Lemma 4.7, we can fix a $B\in K(P\times Q)$ such that
$$B=\bigcap_{i\in I}e_P(A_i)=\bigcap_{i\in I}(A_i\times\up_Qq_0)=(\bigcap_{i\in I}A_i)\times\up_Qq_0.$$
Consequently, $\pi_P(B)=\bigcap\limits_{i\in I}A_i$ is a nonempty compact saturated subset of $\Sigma P$ as $\pi_P$ is Scott continuous. By Lemma 4.6 and Lemma 4.7 again, $\bigcap\limits_{i\in I}A_i=\bigvee\limits_{i\in I}A_i$. Therefore $K(P)$ is a dcpo. Similarly, $K(Q)$ is also a dcpo.

Define $r_P:K(P\times Q)\longrightarrow K(P)$ by $r_P(B)=\pi_P(B)$.
Then $\pi_P(B)$ is a nonempty compact saturated subset of $\Sigma P$. Thus $r_P$ is well defined.
One can check that $e_{P}$ is Scott continuous and $r_P\circ e_P=\mathrm{id}_{K(P)}$.

Let $\{B_j\mid j\in J\}$ be a nonempty filtered family in $K(P\times Q)$. Clearly, $\pi_P(\bigcap\limits_{j\in J}B_j)\subseteq\bigcap\limits_{j\in J}\pi_P(B_j)$. Conversely, choose an arbitrary $p\in\bigcap\limits_{j\in J}\pi_P(B_j)$ and put
$$F_j=\{q\in Q\mid(p,q)\in B_j\}.$$
Each $F_j$ is nonempty and saturated in $Q$. In addition,
$$F_j=\pi_Q(B_j\cap(\dn_P p\times Q)),$$ because $B_{j}$ is upper.

Since $\downarrow_Pp$ is Scott closed, $\dn_Pp\times Q=\pi_P^{-1}(\dn_Pp)$
is Scott closed in $P\times Q$. So we have
$B_j\cap(\dn_Pp\times Q)$ is compact
and hence $F_j\in K(Q)$. The family $\{F_\lambda:\lambda\in\Lambda\}$ is filtered as the family $\{B_j\mid j\in J\}$ is filtered.
As $K(Q)$ is a dcpo, by Lemma 4.6, we have
$\bigcap_{j\in J}F_j$ is a nonempty compact saturated subset of $\Sigma Q$.
Choose $$q\in\bigcap_{j\in J}F_j.$$ This leads to
$(p,q)\in\bigcap_{j\in J}B_j$ and hence $p\in\pi_P(\bigcap_{j\in J}B_j)$.
As a result, $\pi_P(\bigcap\limits_{j\in J}B_j)=\bigcap\limits_{j\in J}\pi_P(B_j)$. So $r_P$ is Scott continuous. Consequently, $\Sigma K(P)$ is a retract of $\Sigma K(P\times Q)$. By Lemma 3.6, $\Sigma K(P)$ is well-filtered.\end{proof}

Naturally, it arises the following problem:

\begin{problem}
{\rm Let $P$ and $Q$ be dcpos such that $\Sigma K(P),\Sigma K(Q)$ are well-filtered and $K(P\times Q)$ is a dcpo. Is $\Sigma K(P\times Q)$ well-filtered?}
 \end{problem}

In the following, we will give a partial answer to the above problem (see Theorem 4.17).

\begin{lemma}
{\rm (\cite{GG03}) A closed subspace of a well-filtered space is well-filtered.}
\end{lemma}

\begin{lemma}{\rm (\cite{HMOOP281})
 Let $P$ be a dcpo and $A$ a nonempty Scott closed subset of $P$. Then $A$ is a subdcpo and its Scott topology is the subspace topology inherited from $\Sigma P$.}
\end{lemma}

\begin{lemma}{\rm \cite{XSXZ20} Let $\mathscr{K}$ be a nonempty family consisting of nonempty compact subsets of a topological space $X$ and $F\subseteq X$ a nonempty closed subset such that $F\in M(\mathscr{K})$. Then there is a closed set
$A$ contained in $F$ such that $A$ is a minimal element of $M(\mathscr{K})$. In particular, if $\mathscr{K}\subseteq K(X)$ is filtered, then $A$ is an irreducible closed subset of $X$.}
\end{lemma}

\begin{lemma}{\rm (\cite{XSXZ20})
Let $f:X\longrightarrow Y$ be a continuous map between $T_0$ spaces $X$ and $Y$. If
$A\in {\rm RD}(X)$, then $cl_Y (f(A))\in RD(Y)$.}
\end{lemma}

\begin{lemma}
{\rm Let $P$ be a dcpo and $\{q_i\mid i\in I\}$ be a directed set consisting of Scott continuous self-maps on $P$ such that $\bigvee\limits_{i\in I}q_{i}=id_{P}$.
 If there exist a Scott well-filtered dcpo $W_{i}$ and two Scott continuous mappings $f_{i}:P\longrightarrow W_{i}$ and $g_{i}:W_{i}\longrightarrow P$ such that $q_i=g_i\circ f_i$ for each $i$,
then $\Sigma P$ is well-filtered.}
\end{lemma}

\begin{proof}
Let $A\subseteq P$ be a closed Rudin set. By Lemma 4.13, $cl_{\Sigma W_i} (f_i(A))$ is a closed Rudin set in $\Sigma W_{i}$.  As $\Sigma W_i$ is well-filtered, by Theorem 3.3, there is $w_i\in W_i$ with
$$
cl_{\Sigma W_i} (f_i(A))=\dn_{W_i}w_i.
$$
Take $a_i=g_i(w_i)$. As $g_{i}$ is an order preserving, we have $$q_i(A)=g_{i}(f_{i}(A))\subseteq  g_{i}(\dn_{W_i}w_i)\subseteq\dn_P a_i.$$ Since $g_i$ is Scott continuous, we know
$$a_i\in g_i(cl_{\Sigma W_i}f_i(A))\subseteq cl_{\Sigma P}(q_i(A)).$$
As $\downarrow_P a_i$ is Scott closed, we deduce that  $cl_{\Sigma P}(q_i(A))=\downarrow_P a_i$.
It follows from $q_i\leq id_P$ and the Scott closedness of $A$ that $q_i(A)\subseteq A$. This implies $a_i\in A$.

Choose any $q_i,q_j$. We can select a $q_{k}$ with $q_{i}, q_{j}\leq q_{k}$. Then we have $q_i(A),q_j(A)\subseteq\dn_P a_k$. As a consequence, $\dn_{P}a_i=cl_{\Sigma P}(q_i(A))\subseteq\dn a_{k}$ and $\dn_{P}a_j=cl_{\Sigma P}(q_j(A))\subseteq\dn a_{k}$. This yields that $\{a_i\mid i\in I\}$ is directed. Take $a=\bigvee\limits_{i\in I} a_i$. Then $a\in A$ since $A$ is Scott closed. For every $x\in A$,
$$
x=\bigvee\limits_{i\in I}q_i(x)\leq\bigvee\limits_{i\in I} a_i=a.
$$
Thus $A=\dn_P a$ because $A$ is lower. By Theorem 3.3, $\Sigma P$ is well-filtered.
\end{proof}

\begin{lemma}
{\rm Let $M$ be a dcpo such that $\Sigma K(M)$ is well-filtered. If $F$ is a finite poset, then $\Sigma K(M\times F)$ is well-filtered.}
\end{lemma}

\begin{proof} Set
$$
S=K(M),\text{ and } T=S\cup\{\varnothing\},
$$
both ordered by the reverse inclusion. As $\Sigma K(M)$ is well-filtered, $S$ is a dcpo. Hence $T$ is also a dcpo with a greatest element $\top_T=\varnothing$ and further $\top_T$ is a compact element. It is not difficult to check that $$K(T)=\{Q\cup\{\top_T\}\mid Q\in K(S)\cup\{\emptyset\}\}$$ and $$\sigma(T)=\{A\cup\{\top_T\}\mid A\in\sigma(S)\}\cup\{\emptyset\}.$$ It follows from the well-filteredness of $\Sigma S$ that $\Sigma T$ is well-filtered. We enumerate $F$ as $\{a_{1},a_{2},\cdot\cdot\cdot,a_{n}\}$. Let
$$
B=T^n,
\top_{B}=(\top_{T},\top_{T},\cdot\cdot\cdot,\top_{T})\in B\text{ and }B_0=B\setminus\{\top_{B}\}.$$
By Theorem 3.5, $\Sigma B$ is well-filtered. The greatest element $\top_{B}$ is compact in $B$ because $\top_{T}$ is compact in $T$. So $\{\top_{B}\}$ is Scott open in $B$ and hence $B_{0}$ is Scott closed. By Lemma 4.10 and Lemma 4.11, $\Sigma B_0$ is well-filtered. For $D=(D_{1},D_{2},\cdot\cdot\cdot,D_{n})\in B$ and $k\in\{1,\cdot\cdot\cdot,n\}$ define
$$p_{D}(k)=\bigcup\{D_{m}\mid m\in\{1,\cdot\cdot\cdot,n\}, a_{m}\leq_F a_{k}\}.$$
Then each $p_{D}(k)\in T$. The mapping $p:B\longrightarrow B$ is defined by $$\forall D\in B,\ p(D)=(p_{D}(1),p_{D}(2),\cdot\cdot\cdot,p_{D}(n)).$$ Then we can check that $p$ is an order preserving and $p(\top_{B})=\top_{B}$.

Let $\{D_{i}\mid i\in I\}\subseteq B$ be a directed subset. For convenience, put $D_{i}=(D_{i,1},D_{i,2},\cdot\cdot\cdot D_{i,n})$ for each $i\in I$. Then by Lemma 4.6,

$$
\begin{aligned}
&\quad p(\bigvee\limits_{i\in I} D_{i})\\
&\quad=p(\bigvee\limits_{i\in I} D_{i,1},\bigvee\limits_{i\in I} D_{i,2},\cdot\cdot\cdot,\bigvee\limits_{i\in I} D_{i,n})\\
&\quad=p(\bigcap\limits_{i\in I} D_{i,1},\bigcap\limits_{i\in I} D_{i,2},\cdot\cdot\cdot,\bigcap\limits_{i\in I} D_{i,n})\\
&\quad=(\bigcup\limits_{a_{m}\leq_{F}a_{1}}\bigcap\limits_{i\in I} D_{i,m},\bigcup\limits_{a_{m}\leq_{F}a_{2}}\bigcap\limits_{i\in I} D_{i,m},\cdot\cdot\cdot,\bigcup\limits_{a_{m}\leq_{F}a_{n}}\bigcap\limits_{i\in I} D_{i,m})
\end{aligned}
$$
and
$$
\begin{aligned}
&\quad \bigvee\limits_{i\in I}p(D_{i})\\
&\quad=\bigvee\limits_{i\in I}(p_{D_{i}}(1),p_{D_{i}}(2),\cdot\cdot\cdot,p_{D_{i}}(n))\\
&\quad=(\bigcap\limits_{i\in I} p_{D_{i}}(1),\bigcap\limits_{i\in I} p_{D_{i}}(2),\cdot\cdot\cdot,\bigcap\limits_{i\in I} p_{D_{i}}(n))\\
&\quad=(\bigcap\limits_{i\in I}\bigcup\limits_{a_{m}\leq_{F}a_{1}}D_{i,m},\bigcap\limits_{i\in I} \bigcup\limits_{a_{m}\leq_{F}a_{2}}D_{i,m},\cdot\cdot\cdot,\bigcap\limits_{i\in I} \bigcup\limits_{a_{m}\leq_{F}a_{n}}D_{i,m}).
\end{aligned}
$$
We claim that for each $1\leq k\leq n$, $$\bigcup\limits_{a_{m}\leq_{F}a_{k}}\bigcap\limits_{i\in I} D_{i,m}=\bigcap\limits_{i\in I} \bigcup\limits_{a_{m}\leq_{F}a_{k}}D_{i,m}.$$

Choose any $r\in\bigcup\limits_{a_{m}\leq_{F}a_{k}}\bigcap\limits_{i\in I} D_{i,m}$. Select a $m_{0}\in\{1,2,\cdot\cdot\cdot,n\}$ with $a_{m_{0}}\leq_{F} a_{k}$ such that $r\in\bigcap\limits_{i\in I} D_{i,m_{0}}$. As for each $i\in I$, $D_{i,m_{0}}$ belongs to $\bigcup\limits_{a_{m}\leq_{F}a_{k}}D_{i,m}$, we have $r\in\bigcap\limits_{i\in I} \bigcup\limits_{a_{m}\leq_{F}a_{k}}D_{i,m}$. This yields $\bigcup\limits_{a_{m}\leq_{F}a_{k}}\bigcap\limits_{i\in I} D_{i,m}\subseteq\bigcap\limits_{i\in I} \bigcup\limits_{a_{m}\leq_{F}a_{k}}D_{i,m}$.

Conversely, let $s\in\bigcap\limits_{i\in I} \bigcup\limits_{a_{m}\leq_{F}a_{k}}D_{i,m}$. Then for each $i\in I$, we can fix a $a_{m_{i}}\leq a_{k}$ with $s\in D_{i,m_{i}}$. Assume that $s\not\in\bigcup\limits_{a_{m}\leq_{F}a_{k}}\bigcap\limits_{i\in I} D_{i,m}$. Put $$X=\{p\mid 1\leq p\leq n, a_{p}\leq_{F} a_{k}\}.$$ It follows that for each $p\in X$, $s\not\in\bigcap\limits_{i\in I} D_{i,p}$. Select a $i_{p}\in I$ with $s\not\in D_{i_{p},p}$. As $\{D_{i}\mid i\in I\}$ is directed and $X$ is finite, there is a $j\in I$ such that $D_{j}\subseteq \bigcap\limits_{p\in X}D_{i_{p}}$. For $m_{j}\in X$, we have $s\in D_{j,m_{j}}$ and $s\not\in D_{i_{m_{j}},m_{j}}$. It follows from $D_{j}\subseteq \bigcap\limits_{p\in X}D_{i_{p}}$ that $D_{j,m_{j}}\subseteq \bigcap\limits_{p\in X}D_{i_{p},m_{j}}$. As a result, $s\in\bigcap\limits_{p\in X}D_{i_{p},m_{j}}$, a contradiction. So we know $s\in\bigcup\limits_{a_{m}\leq_{F}a_{k}}\bigcap\limits_{i\in I} D_{i,m}$. Thus $\bigcap\limits_{i\in I} \bigcup\limits_{a_{m}\leq_{F}a_{k}}D_{i,m}\subseteq\bigcup\limits_{a_{m}\leq_{F}a_{k}}\bigcap\limits_{i\in I} D_{i,m}$.

Hence $p$ is Scott continuous.

For each $E=(E_{1},E_{2},\cdot\cdot\cdot,E_{n})\in B$ and each $1\leq k\leq n$,
$$p_{p(E)}(k)=\bigcup\{p_{E}(m)\mid a_{m}\leq_{F}a_{k}\}=\bigcup\limits_{a_{m}\leq_{F}a_{k}}\bigcup\limits_{a_{q}\leq_{F}a_{m}}E_{q}=\bigcup\limits_{a_{q}\leq_{F}a_{k}}    E_{q}=p_{E}(k).$$ This means that $p$ is idempotent.

Set
$$\mathscr A=\{G\in B_{0}\mid \text{ if }a_{i}\leq_F a_{j},\text{ then } G_{i}\subseteq G_{j}\}.$$

Let $G=(G_1,G_{2},\ldots,G_n)\in B_0$ with $p(G)=G$. Suppose $a_j\leq_F a_k$. Then we have
$$
G_j\subseteq
\bigcup_{\substack{a_m\leq_F a_k}}G_m
=p_G(k)=G_k.
$$
So $G\in\mathscr A$.

Suppose $G\in\mathscr A$. Then for each $1\leq k\leq n$, $G_m\subseteq G_k$ whenever $a_m\leq_F a_k$. This yields that
$$
p_G(k)=
\bigcup\limits_{a_m\leq_F a_k}G_m
=G_k.
$$
Consequently, $p(G)=G$.

We conclude that $\mathscr A=\{G\in B_0\mid p(G)=G\}$. Since $p$ is idempotent,
we deduce that $\{p(G)\mid G\in B_{0}\}=\mathscr A$.

Suppose $\mathscr D$ is  directed
in $\mathscr A$ and set $H=\bigvee_{B}\mathscr D$. Since $\top_{B}$ is compact in $B$, $H=\bigvee_{B}\mathscr D=\bigvee_{B_{0}}\mathscr D\in B_{0}$. By the Scott continuity of $p$, we obtain
$$
p(H)
=p\left(\bigvee_{B_{0}}\mathscr D\right)
=\bigvee_{B_{0}}\{p(D)\mid D\in\mathscr D\}
=\bigvee_{B_{0}}\mathscr D
=H.
$$
Hence $H\in\mathscr A$. As a result, $\mathscr A$ is a subdcpo
of $B_0$.

Let $i:\mathscr A\to B_0$ be the inclusion mapping. Then $i$
is Scott continuous. Set $q=p|B_0$ as the mapping $p|B_0:B_0\to\mathscr A$. Also, $q$ is
Scott continuous. Moreover, $q\circ i=id_{\mathscr A}$.

Since $\Sigma B_0$ is well-filtered, by Lemma 4.14,
$\Sigma\mathscr A$ is well-filtered.

Finally, we verify that $K(M\times F)$ and $\mathscr A$ are order isomorphic.

For $Q\in K(M\times F)$ and
$1\leq k\leq n$, put
$$
Q_k=\{x\in M\mid(x,a_k)\in Q\}.
$$
Since $Q$ is upper, we have $Q_k=\pi_M(Q\cap(M\times\dn_F a_k))$. Clearly, the projections $\pi_{M}$ and $\pi_{F}$ are Scott continuous. As $\dn_F a_k$ is closed in $\Sigma F$,
we get
$$
M\times\dn_F a_k
=\pi_F^{-1}(\dn_F a_k)
$$
is Scott closed in $M\times F$. Thus
$Q\cap(M\times\downarrow_F a_k)$ is compact in $\Sigma(M\times F)$. The Scott continuity of $\pi_{M}$ gives $Q_k$ is compact in $\Sigma M$.
Moreover, $Q_k$ is upper in $M$. So $Q_k\in T$.

If $a_j\leq_F a_k$, then we have $Q_j\subseteq Q_k$ as $Q$ is upper.
Since $Q$ is nonempty, at least one of $Q_1,\ldots,Q_n$ is nonempty.
So we can define a mapping $\Phi:K(M\times F)\longrightarrow\mathscr A$ by
$$\forall\ Q\in K(M\times F),\ \Phi(Q)=(Q_1,\ldots,Q_n).$$

For each $G=(G_1,\ldots,G_n)\in\mathscr A$, put
$$
C_G=\{(x,a_k)\mid 1\leq k\leq n,\ x\in G_k\}.
$$

Clearly, $C_G\subseteq\bigcup\limits_{k=1}\limits^{n}
\bigl(G_k\times\up_F a_k\bigr)$. Conversely, suppose
$(x,a_j)\in G_k\times\up_F a_k$. Then $x\in G_k$
and $a_k\leq_F a_j$.  $G\in\mathscr A$
gives $G_k\subseteq G_j$. So $x\in G_j$ and therefore
$(x,a_j)\in C_G$. This shows that $$
C_G=\bigcup\limits_{k=1}\limits^{n}
(G_k\times\up_F a_k).
$$

Since each $G_k$ is an upper set in $M$, each
$G_k\times\uparrow_F a_k$ is also upper in $M\times F$.
So $C_G$ is upper. Besides, $C_G$ is nonempty because
$G\in B_0$.

For each $1\leq k\leq n$, define a mapping
$j_k:M\longrightarrow M\times F$ by $j_k(x)=(x,a_k)$.
Then $j_{k}$ preserves the directed superma and therefore is Scott continuous.
The compactness of $G_k$ gives 
$$
j_k(G_k)=G_k\times\{a_k\}
$$
is compact in $\Sigma(M\times F)$. So its saturation 
$$
\up_{M\times F}(G_k\times\{a_k\})
=G_k\times\up_F a_k,
$$
is also compact in $\Sigma(M\times F)$. As a consequence, $C_G$ is a finite union of compact sets
and is therefore compact. In other words, $C_G\in K(M\times F)$ and hence the mapping
$$\Psi:\mathscr A\longrightarrow K(M\times F),\ \Psi(G)=C_G$$ is well defined.

By construction, for every $G\in\mathscr A$ and every
$1\leq k\leq n$,
$$
\{x\in M\mid(x,a_k)\in C_G\}=G_k.
$$
Hence $\Phi(\Psi(G))=G$. Conversely, for each 
$C\in K(M\times F)$,
$$
C_{\Phi(C)}
=\{(x,a_k)\mid 1\leq k\leq n,\ x\in C_k\}
=C.
$$
Therefore $\Phi$ and $\Psi$ are mutually inverse bijections.

Furthermore, for every $C,E\in K(M\times F)$,
$$
C\supseteq E
\quad\Longleftrightarrow\quad
C_k\supseteq E_k
\quad\text{for every }1\leq k\leq n.
$$
Consequently, $\Phi$ is an order isomorphism and hence
$$
K(M\times F)\cong\mathscr A.
$$
Therefore, $\Sigma K(M\times F)$ is well-filtered.\end{proof}

\begin{lemma}
{\rm Let $A,B,R$ be dcpos such that $\Sigma B$ is well-filtered and $K(A\times B)$,$K(A\times R)$ are dcpos. Suppose $f:B\to R$ is Scott continuous. Define a mapping
$\Phi_f:K(A\times B)\longrightarrow K(A\times R)$ by
$$\Phi_f(C)=\up_{A\times R}\{(a,f(b))\mid (a,b)\in C\}.$$
Then $\Phi_f$ is Scott continuous.}
\end{lemma}

\begin{proof}
Define a mapping $h=\operatorname{id}_A\times f: A\times B\longrightarrow A\times R$ by
$h(a,b)=(a,f(b))$.

Let $\{(a_i,b_i):i\in I\}$ be a directed subset of $A\times B$. Then by the Scott continuity of $f$, we have
$$
\begin{aligned}
h\left(\bigvee_{i\in I}(a_i,b_i)\right)
&=h\left(\bigvee_{i\in I}a_i,\bigvee_{i\in I}b_i\right)\\
&=\left(\bigvee_{i\in I}a_i,
f\left(\bigvee_{i\in I}b_i\right)\right)\\
&=\left(\bigvee_{i\in I}a_i,\bigvee_{i\in I}f(b_i)\right)\\
&=\bigvee_{i\in I}h(a_i,b_i).
\end{aligned}
$$
Thus $h$ is Scott continuous.

Let $C\in K(A\times B)$. Then $\Phi_f(C)=\up_{A\times R}h(C)$
is a nonempty compact saturated subset of $\Sigma(A\times R)$.
Hence $\Phi_f$ is well defined.

Moreover, if $C,D\in K(A\times B)$ with $C\supseteq D$, then
$$
\Phi_f(C)=\up_{A\times R}h(C)
\supseteq
\up_{A\times R}h(D)=\Phi_f(D).
$$
Consequently, $\Phi_f$ is an order preserving. For $C\in K(A\times B)$ and $a\in A$, set
$$
C_a=\{b\in B\mid (a,b)\in C\}.
$$
We claim that $C_a=\pi_B(C\cap(\dn_A a\times B))$.

Clearly $C_a\subseteq\pi_B(C\cap(\dn_A a\times B))$. Conversely, if $b$ belongs to $\pi_B(C\cap(\dn_A a\times B))$,
then there exists $a'\leq_A a$ such that $(a',b)\in C$.
As $C$ is upper and $(a',b)\leq(a,b)$, we obtain
$(a,b)\in C$ and therefore $b\in C_a$. As a result, $\pi_B(C\cap(\dn_A a\times B))\subseteq C_{a}$. This proves that $C_a=\pi_B(C\cap(\dn_A a\times B))$.

Since $\dn_A a$ is Scott closed and $\pi_{A}$ is Scott continuous, we get
$$
\dn_A a\times B
=\pi_A^{-1}(\dn_A a)
$$
is Scott closed in $A\times B$. Hence
$C\cap(\dn_A a\times B)$ is also compact in $\Sigma(A\times B)$. It follows from the Scott continuity of $\pi_{B}$ that $C_a=\pi_B(C\cap(\dn_A a\times B))$ is compact in
$\Sigma B$. Furthermore, $C_a$ is an upper set in $B$ because $C$ is upper in $A\times B$. 

For every $a\in A$ and $r\in R$, one can check that 
$$
(a,r)\in\Phi_f(C)
\quad\Longleftrightarrow\quad
C_a\cap f^{-1}(\dn_R r)\neq\varnothing.
$$
Indeed, suppose that $(a,r)\in\Phi_f(C)$. Then there
exists $(a',b)\in C$ such that
$$
a'\leq_A a
\text{ and }
f(b)\leq_R r.
$$
Since $C$ is upper, $(a,b)\in C$ and hence
$b\in C_a\cap f^{-1}(\dn_R r)$. Conversely, if $b\in C_a\cap f^{-1}(\dn_R r)$, then
$(a,b)\in C$ and $f(b)\leq_R r$. Hence
$$
(a,f(b))\in h(C)
\text{ and }
(a,f(b))\leq(a,r).
$$
This implies $(a,r)\in\Phi_f(C)$. 

Finally, we prove that $\Phi_f$ preserves directed suprema.

Let $\{C_\lambda\mid\lambda\in\Lambda\}$ be a directed subset
of $K(A\times B)$. As $K(A\times B)$ is a dcpo,
by Lemma 4.6, we have that
$$
C=\bigcap\limits_{\lambda\in\Lambda}C_\lambda
=\bigvee\limits_{\lambda\in\Lambda}\limits^{K(A\times B)}C_\lambda
$$
belongs to $K(A\times B)$.

Since $C\subseteq C_\lambda$ for every $\lambda\in\Lambda$,
we have $\Phi_f(C)\subseteq\Phi_f(C_\lambda)$ and hence
$\Phi_f(C)\subseteq\bigcap_{\lambda\in\Lambda}\Phi_f(C_\lambda)$.
Conversely, take any $(a,r)\in\bigcap_{\lambda\in\Lambda}\Phi_f(C_\lambda)$
and let $H=f^{-1}(\dn_R r)$.
The Scott continuity of $f$ gives that $H$ is Scott closed in $B$.

By the characterization established above, for every
$\lambda\in\Lambda$, $(C_\lambda)_a\cap H\neq\varnothing$.
In particular, every $(C_\lambda)_a$ is nonempty and hence for every
$\lambda\in\Lambda$, $(C_\lambda)_a\in K(B)$. 
Moreover, $\{(C_\lambda)_a\mid \lambda\in\Lambda\}$ forms a filtered family. Assume that $H\cap(\bigcap_{\lambda\in\Lambda}(C_\lambda)_a)=\varnothing$.
Equivalently, $\bigcap_{\lambda\in\Lambda}(C_\lambda)_a\subseteq B\setminus H$. 
Since $\Sigma B$ is
well-filtered, there exists $\lambda_0\in\Lambda$ satsfying 
$(C_{\lambda_0})_a\subseteq B\setminus H$.
This contradicts 
$(C_{\lambda_0})_a\cap H\neq\varnothing$. So we obtain 
$H\cap(\bigcap\limits_{\lambda\in\Lambda}(C_\lambda)_a)\neq\varnothing$.

In addition, we can check that 
$$
\begin{aligned}
C_a
&=\{b\in B:(a,b)\in C\}\\
&=\left\{b\in B:
(a,b)\in C_\lambda
\text{ for every }\lambda\in\Lambda\right\}\\
&=\bigcap_{\lambda\in\Lambda}(C_\lambda)_a.
\end{aligned}
$$
This yields that $C_a\cap f^{-1}(\downarrow_R r)=C_a\cap H\neq\varnothing$.
By the characterization above again, we obtain
$(a,r)\in\Phi_f(C)$. Consequently, we know
$$
\Phi_f\left(\bigcap_{\lambda\in\Lambda}C_\lambda\right)
=
\bigcap_{\lambda\in\Lambda}\Phi_f(C_\lambda).
$$

Clearly, the subfamily
$\{\Phi_f(C_\lambda)\mid\lambda\in\Lambda\}$ is directed in
$K(A\times R)$. By Lemma 4.6 and the above equality, we have
$$
\Phi_f(\bigvee_{\lambda\in\Lambda}C_\lambda)
=\Phi_f(\bigcap_{\lambda\in\Lambda}C_\lambda)
=\bigcap_{\lambda\in\Lambda}\Phi_f(C_\lambda)
=\bigvee_{\lambda\in\Lambda}\Phi_f(C_\lambda).
$$
Therefore $\Phi_f$ is Scott continuous.
\end{proof}

\begin{definition}{\rm (\cite{E2}) An $RB$-domain $L$ is a dcpo that carries a directed set consisting of Scott continuous mappings $f_{i}:L\longrightarrow L$, $i\in I$ such that each image set $f_{i}(L)$ is finite and $\bigvee\limits_{i\in I}f_{i}=id_{L}$.}
\end{definition}

\begin{lemma}{\rm(\cite{XWX23}) Let $X$ be a well-filtered $T_{0}$ space. Then $\Sigma K(X)$ is well-filtered.}

\end{lemma}

\begin{remark} {\rm Every $RB$-domain $L$ is a continuous dcpo. So $\Sigma L$ is sober and hence well-filtered. By Lemma 4.18, $\Sigma K(L)$ is well-filtered.
}
\end{remark}

\begin{theorem}
{\rm Let $M$ and $N$ be dcpos such that $\Sigma K(M)$ is well-filtered and $K(M\times N)$ is a dcpo. If There exist a directed set $I$ and Scott continuous maps $d_i:N\to N$, $i\in I$, satisfying

(i) each $d_i(N)$ is finite and $d_i\leq id $;

(ii) $d_i\leq d_j$ if $i\leq j$;

(iii) $\bigvee\limits_{i\in I}d_i=id$.

Then $\Sigma K(M\times N)$ is well-filtered.}
\end{theorem}

\begin{proof}
By Definition 4.17, $N$ is an $RB$-domain. Hence by Remark 4.19, $\Sigma N$ is well-filtered. For each $i\in I$, set $N_{i}=d_{i}(N)$. Let $\bar{d_{i}}$ be the corestriction of the mapping $d_{i}$ and $j_{i}:N_{i}\longrightarrow N$ the inclusion mapping. Then for each $i$, both $\bar d_i$ and $j_i$ are Scott continuous with $d_{i}=j_{i}\circ\bar d_i$.

Set $Q=K(M\times N)$ and for each $i\in I$,  $Q_i=K(M\times N_i)$. By Lemma 4.15,
$Q_i$ is a dcpo and $\Sigma Q_i$ is well-filtered. Define two mappings
$$
\begin{aligned}
p_i:Q&\longrightarrow Q_i,
&p_i(C)&=\up_{M\times N_i}
(id_M\times\bar d_i)(C),\\
r_i:Q_i&\longrightarrow Q,
&r_i(E)&=\up_{M\times N}
(id_M\times j_i)(E).
\end{aligned}
$$
By Lemma 4.16, each $p_i$ and $r_i$ are Scott continuous. 
Let $q_i=r_i\circ p_i$. Then every $q_i$ is also Scott continuous.

Choose any $C\in Q$ and $(x,y)\in C$. It follows from $d_i(y)\leq_N y$ that
$(x,d_i(y))\leq(x,y)$. So we have $C\subseteq q_i(C)$.
If $i\leq j$, then $d_i(y)\leq_N d_j(y)$ for every $y\in N$. So $q_j(C)\subseteq q_i(C)$ and hence
$$
q_i\leq q_j\leq id_Q\qquad(i\leq j).
$$
Since $I$ is directed, $\{q_i\mid i\in I\}$ is directed under the
pointwise order.

Fix $C\in Q$. We show that $\bigcap\limits_{i\in I}q_i(C)=C$.

Clearly, $C\subseteq\bigcap\limits_{i\in I}q_i(C)$. Suppose $W$ is Scott open subset of $M\times N$ containing $C$. For each $i\in I$, take
$$
W_i=(id_M\times d_i)^{-1}(W).
$$
Each $W_i$ is a Scott open subset of $M\times N$ since $id_M\times d_i$ is Scott continuous. $\forall\ i,j\in I$ with $i<j$, choose any $(a,b)\in W_{i}$ i.e. $(a,d_{i}(b))\in W$. As $d_{i}\leq d_{j}$, $(a,d_{i}(b))\leq(a,d_{j}(b))$. It follows from $W\in\sigma(M\times N)$ that $(a,d_{j}(b))\in W$.  So $(a,b)\in W_{j}$ and therefore $W_{i}\subseteq W_{j}$. 
Take any $(x,y)\in C$, the condition ${\rm(iii)}$ gives
$$
(x,y)=\sup\limits_{i\in I}(x,d_i(y)).
$$
Since $(x,y)\in C\subseteq W$ and  $W\in\sigma(M\times N)$, there is 
$(x,d_i(y))\in W$, or equivalently $(x,y)\in W_i$.
This means that $$C\subseteq\bigcup\limits_{i\in I}W_i.$$ By the compactness of $C$, there exists $i_{1},i_{2},\cdot\cdot\cdot,i_{n}\in I$ such that $C\subseteq\bigcup\limits_{k=1}\limits^{n}W_{i_{k}}$.
The directedness of $I$ provides $i_0$ being an upper bound of $\{i_{1},i_{2},\cdot\cdot\cdot,i_{n}\}$. As a consequence, $C\subseteq W_{i_0}$.
In other words, $(id_M\times d_{i_0})(C)\subseteq W$. As $W$ is upper, we have
$$\bigcap\limits_{i\in I}q_i(C)\subseteq q_{i_0}(C)=r_{i_{0}}(p_{i_0}(C))=\up_{M\times N}
(id_M\times j_{i_{0}})((\up_{M\times N_{i_{0}}}
id_M\times\bar d_{i_{0}})(C))\subseteq W.$$

Since $C$ is saturated, we conclude that $\bigcap\limits_{i\in I}q_i(C)\subseteq C$. This shows that $\bigcap\limits_{i\in I}q_i(C)=C$. By Lemma 4.6, we have $\bigvee\limits_{i\in I}q_{i}=id_{Q}$. By Lemma 4.14, $\Sigma Q=\Sigma K(M\times N)$ is well-filtered.
\end{proof}

\begin{remark}{\rm In Theorem 4.20, the dcpo $N$ is actually a ${\rm RB}$-domain. Consequently, if $M$ is a dcpo such that $\Sigma K(M)$ is well-filtered and $N$ an ${\rm RB}$-domain, then $\Sigma K(M\times N)$ is well-filtered whenever $K(M\times N)$ is a dcpo.}
\end{remark} 

\section{An intrinsic characterization of WD sets by open filters}

In this section, we will present a characterization of WD sets in a $T_{0}$ space.
\begin{definition}{\rm(\cite{XSXZ20})
A subset $A$ of a $T_0$ space $X$ is called a \emph{well-filtered determined set}, or simply a
\emph{WD set} if for every well-filtered space $Y$ and every continuous mapping $f:X\longrightarrow Y$, there is a unique $y\in Y$ such that $cl_Y (f(A))= cl_Y(\{y\})$.}
\end{definition}

\begin{definition}
{\rm Let $X$ be a $T_{0}$ space. A \emph{proper filter}
$\mathscr{F}\subseteq \mathcal{O}(X)$ is a filtered upper set with $X\in \mathscr{F}$ and $\varnothing\notin \mathscr{F}$. We denote by $\mathsf{SF}(X)$ the set of all proper Scott open filters on $\mathcal{O}(X)$.}

\end{definition}

For $A\subseteq X$, define
$$
 p_A=\{U\in\mathcal O(X)\mid U\cap A\ne\varnothing\}.
$$
For $x\in X$, we write $p_{\{x\}}$ as $p_x$.

\begin{lemma}
{\rm Let $X$ be a $T_{0}$ space and $A\subseteq X$. Then the following statements hold:

(1) $p_A=p_{cl_X (A)}$;

(2) If $A$ is nonempty, then $p_A$ is a nonempty Scott open set in $\mathcal{O}(X)$.

(3) $p_A\in\mathsf{SF}(X)$ if and only if $cl_X (A)$ is an irreducible subset of $X$.}
\end{lemma}

\begin{proof}
(1) The first assertion holds since an open set $U\in\mathcal{O}(X)$ meets $A$ if and only if $U\cap cl_X (A)$ is nonempty.

(2) Clearly, $p_A$ is upper and nonempty in $\mathcal{O}(X)$. Let $\{U_i\mid i\in I\}$ be a subset of $\mathcal{O}(X)$ with $\bigcup\limits_{i\in I}U_{i}\in p_{A}$. Then there is $U_i$ such that $U_{i}\cap A\neq\emptyset$. This proves the Scott openness of $p_{A}$.

(3) Let $cl_{X}(A)$ be an irreducible subset of $X$. Then $A$ is also irreducible. Choose any $U,V\in p_{A}$. Then $U\cap A\neq\emptyset$ and $V\cap A\neq\emptyset$. The irreducibility of $A$ gives $U\cap V\in p_{A}$. This yields that $p_{A}$ is filtered.  As $A$ is nonempty, by the second assertion, $p_{A}$ is a proper Scott open set. So $p_{A}\in\mathsf{SF}(X)$.

Conversely, suppose $M,N$ are open in $X$ such that $M\cap A\neq\emptyset$ and $N\cap A\neq\emptyset$. In other words, $M,N\in p_{A}$. As $p_{A}$ is filtered and upper, we deduce that $M\cap N\in p_{A}$. So $A$ is irreducible in $X$. The irreducibility of $cl_{X}(A)$ is obtained immediately.
\end{proof}

\begin{definition}
{\rm Let $X$ be a $T_{0}$ space. The collection $\mathscr H_X$ is the smallest subfamily of $\mathsf{SF}(X)$ satisfying the following conditions:

(1) $p_x\in\mathscr H_X$ for every $x\in X$;

(2) If $\{\mathscr{F}_i\mid i\in I\}\subseteq\mathscr H_X$ is directed
 under inclusion, then $\bigcup\limits_{i\in I}\mathscr{F}_i\in\mathscr H_X$.

(3) If $\{\mathscr{F}_i:i\in I\}\subseteq\mathscr H_X$ is nonempty and
 $\bigcap\limits_{i\in I}\mathscr{F}_i\in\sigma(\mathcal O(X))$, then
 $\bigcap_{i\in I}\mathscr{F}_i\in\mathscr H_X$.
}
\end{definition}

\begin{remark}
{\rm (1) The family $\mathsf{SF}(X)$ satisfies the three rules in Definition 5.4.

(2) The intersection of all subfamilies satisfying these rules in Definition 5.4 also satisfies
these rules. Therefore the subfamily $\mathscr H_X$ always exists.}
\end{remark}

A subset $K$ of a $T_{0}$ space $X$ is supercompact if each cover $\mathscr{U}\subseteq\mathcal{O}(X)$ of $K$ contains an element $U\in\mathscr{U}$ containing $K$ (\cite{MU13}).

\begin{lemma}
{\rm (\cite{MU13}) A saturated subset $K$ of a $T_{0}$ space $X$ is supercompact if and only if $K=\up x$ for some $x\in X$.}
\end{lemma}

\begin{lemma} {\rm Let $X$ be a $T_{0}$ space and $Y$ a well-filtered space. Fix a continuous mapping $f:X\longrightarrow Y$ and define a subfamily $R_{f}$ of $\mathsf{SF}(X)$ by
$$\mathscr{F}\in R_{f}\Longleftrightarrow\exists K_{\mathscr{F}}\in K(Y)\text{ s.t. }\forall\ V\in\mathcal{O}(Y), (f^{-1}(V)\in\mathscr{F}\text{ iff } K_{\mathscr{F}}\subseteq V).$$
Then $R_{f}$ satisfies the three rules in Definition 5.4.}
\end{lemma}

\begin{proof}

For any $\mathscr{F}\in R_{f}$, we have $$K_{\mathscr{F}}=\bigcap\{V\in\mathcal{O}(Y)\mid K_{\mathscr{F}}\subseteq V\}=\bigcap\{V\in\mathcal{O}(Y)\mid f^{-1}(V)\in\mathscr{F}\}$$ because $K_{\mathscr{F}}$ is saturated. Consequently, $K_{\mathscr{F}_{2}}\subseteq K_{\mathscr{F}_{1}}$ if $\mathscr{F}_{1}\subseteq\mathscr{F}_{2}$.

(1) For each $x\in X$, take $K_{x}=\up_{Y} f(x)$. Then $K_{x}\in K(Y)$. $\forall\ U\in\mathcal{O}(Y)$, we have $$K_{x}\subseteq U\Longleftrightarrow f(x)\in U\Longleftrightarrow x\in f^{-1}(U)\Longleftrightarrow f^{-1}(U)\in p_{x}.$$ Thus each $p_{x}\in R_{f}$.

(2) Let $\{\mathscr{F}_{i}\mid i\in I\}\subseteq R_{f}$ be a directed subfamily. Then for each $i$, there exists a $K_{\mathscr{F}_{i}}\in K(Y)$ such that for $V\in\mathcal{O}(Y)$, $K_{\mathscr{F}_{i}}\subseteq V$ if and only if $f^{-1}(V)\in\mathscr{F}_{i}$. One can check that $\bigcup\limits_{i\in I}\mathscr{F}_{i}$ is a proper Scott open filter in $\mathcal{O}(X)$. Put $K=\bigcap\limits_{i\in I}K_{\mathscr{F}_{i}}$. As $\{\mathscr{F}_{i}\mid i\in I\}\subseteq R_{f}$ is directed, we know that $\{K_{\mathscr{F}_{i}}\mid i\in I\}$ is a filtered family consisting of nonempty compact saturated sets of $Y$. By Lemma 4.7, $K$ is also a nonempty compact saturated subset of $Y$. For each $H\in\mathcal{O}(Y)$, by the well-filteredness of $Y$, we deduce that
$$
 \begin{aligned}
K\subseteq H
&\quad\Longleftrightarrow\exists i\in I, K_{\mathscr{F}_{i}}\subseteq H\\
&\quad\Longleftrightarrow\exists i\in I, f^{-1}(H)\in\mathscr{F}_{i}\\
&\quad\Longleftrightarrow f^{-1}(H)\in\bigcup\limits_{i\in I}\mathscr{F}_{i}.
 \end{aligned}
$$
Therefore, $\bigcup\limits_{i\in I}\mathscr{F}_{i}\in R_{f}$.

(3) Suppose $\{\mathscr{A}_i:i\in I\}\subseteq R_f$ is nonempty
such that $\mathscr{A}=\bigcap\limits_{i\in I}\mathscr{A}_i$ is Scott open in $\mathcal{O}(X)$. Then for every $i$, we can fix a $K_{\mathscr{A}_{i}}\in K(Y)$ satisfying for every $W\in\mathcal{O}(Y)$, $K_{\mathscr{A}_{i}}\subseteq W$ if and only if $f^{-1}(W)\in\mathscr{A}_{i}$. Take $K=\bigcup\limits_{i\in I}K_{\mathscr{A}_{i}}$.
For each $E\in\mathcal{O}(Y)$
$$
 \begin{aligned}
 K\subseteq E
 &\quad\Longleftrightarrow\forall\ i\in I, K_{\mathscr{A}_{i}}\subseteq E\\
 &\quad\Longleftrightarrow\forall\ i\in I, f^{-1}(E)\in\mathscr{A}_i\\
 &\quad\Longleftrightarrow f^{-1}(E)\in \mathscr{A}.
 \end{aligned}
$$
Next, we prove that $K$ is compact in $Y$. Let $\{G_j\mid j\in J\}\subseteq\mathcal{O}(Y)$ be an open cover of $K$. Then we have $f^{-1}(\bigcup\limits_{j\in J}G_j)\in \mathscr{A}$.
The family
$$\mathscr{B}=\{f^{-1}(\bigcup_{j\in \widetilde{J}}G_j)\mid
 \widetilde{J}\subseteq J\text{ is finite}\}
$$
is directed in $\mathcal{O}(X)$ and further $\bigcup\mathscr{B}=f^{-1}(\bigcup\limits_{j\in J}G_j)$. The Scott openness of $\mathscr{A}$ gives that there exists a finite set $\widehat{J}\subseteq J$ such that $f^{-1}(\bigcup\limits_{j\in \widehat{J}}G_j)\in \mathscr{A}$. This yields $K\subseteq\bigcup\limits_{j\in \widehat{J}}G_j$.
This proves the compactness of $K$. Hence we conclude that $\mathscr{A}\in R_f$.\end{proof}

Let $X$ be a $T_{0}$ space. Define $$\mathscr{P}=\{C\mid C\text{ is a nonempty closed subset of }
X\text{ such that }p_{C}\in\mathscr H_X\}.$$
For each $U\in\mathcal{O}(X)$, set $\diamondsuit U=\{C\in \mathscr{P}\mid C\cap U\neq\emptyset\}$  and $\tau=\{\diamondsuit U\mid U\in\mathcal{O}(X)\}$. By Lemma 5.3 (3), each element of $\mathscr{P}$ is irreducible in $X$. This indicates that $\tau$ forms a topology on $\mathscr{P}$. Moreover, each closed subset of $(\mathscr{P},\tau)$ has the form $\Box V=\{C\in\mathscr{P}\mid C\subseteq V\}$ for some closed subset $V$ of $X$. Let $C,D\in\mathscr{P}$ with $C\not\subseteq D$. Then $C\in\diamondsuit(X\setminus D)$ and $D\not\in\diamondsuit(X\setminus D)$. Hence $(\mathscr{P},\tau)$ is $T_{0}$. Besides, for each $\mathscr{A}\subseteq\mathscr{P}$, we can check that $cl_{\tau}(\mathscr{A})=\Box( cl_{X}(\bigcup\mathscr{A}))$.

\begin{lemma} {\rm Let $\mathscr{B}\subseteq\mathscr{P}$ be a nonempty compact subset of $(\mathscr{P},\tau)$. Then $H_{\mathscr{B}}=\bigcap\limits_{B\in\mathscr{B}}p_{B}$ belongs to $\mathsf{SF}(X)$.}

\end{lemma}

\begin{proof} Clearly, $H_{\mathscr{B}}$ is a proper filter of $\mathcal{O}(X)$. We assert that $H_{\mathscr{B}}$ is Scott open in $\mathcal{O}(X)$. Suppose $\{U_{i}\mid i\in I\}\subseteq\mathcal{O}(X)$ is directed with $\bigcup\limits_{i\in I}U_{i}\in H_{\mathscr{B}}$. Then $(\bigcup\limits_{i\in I}U_{i})\cap B$ is nonempty for each $B\in\mathscr{B}$. This means that $\{\diamondsuit U_{i}\mid i\in I\}\subseteq\tau$ covers $\mathcal{B}$. The compactness of $\mathcal{B}$ gives that $$\mathcal{B}\subseteq(\diamondsuit U_{i_{1}})\cup(\diamondsuit U_{i_{2}})\cup\cdot\cdot\cdot\cup(\diamondsuit U_{i_{n}}).$$ As $\{U_{i}\mid i\in I\}$ is directed, we can choose a $U_{i^{\ast}}$ such that $ U_{i_{1}}\cup U_{i_{2}}\cup\cdot\cdot\cdot\cup U_{i_{n}}\subseteq U_{i^{\ast}}$. As a result, $\mathcal{B}\subseteq\diamondsuit U_{i^{\ast}}$. This yields that $U_{i^{\ast}}\in H_{\mathscr{B}}$. The Scott openness of $H_{\mathscr{B}}$ is proved. Therefore, $H_{\mathscr{B}}\in\mathsf{SF}(X)$.
\end{proof}

\begin{lemma}{\rm Let $X$ be a $T_{0}$ space. Then $(\mathscr{P},\tau)$ is well-filtered.}

\end{lemma}

\begin{proof}

Suppose $\{\mathscr{C}_{j}\mid j\in J\}$ is a filtered family consisting of compact saturated subsets and $\mathscr{D}$ is a closed subset in $(\mathscr{P},\tau)$ such that each $\mathscr{C}_{j}\cap\mathscr{D}\neq\emptyset$. By Lemma 4.12, there is a minimal closed subset $\mathscr{E}$ contained in $\mathscr{D}$ that meets each $\mathscr{C}_{j}$. Set $$E=cl_{X}(\bigcup\mathscr{E}),\ \mathscr{E}_{j}=\mathscr{C}_{j}\cap\mathscr{E}\text{ and } H_{j}=\bigcap\limits_{F\in\mathscr{E}_{j}}p_{F}.$$ As $\mathscr{E}$ is closed and $\mathscr{C}_{j}$ is compact, we know that $\mathscr{E}_{j}$ is also compact. By Lemma 5.8, each $H_{j}$ is a proper Scott open filter in $\mathcal{O}(X)$. As $\mathscr{H}_{X}$ satisfies the third rule in Definition 5.4 and for each $F\in\mathscr{E}_{j}$, $p_{F}$ belongs to $\mathscr{H}_{X}$, we know that each $H_{j}$ belongs to $\mathscr{H}_{X}$. As $\{\mathscr{C}_{j}\mid j\in J\}$ is filtered, we deduce that $\{H_{j}\mid j\in J\}$ is directed.

Choose any $C\in\bigcup\limits_{j\in J}H_{j}$. Then we can fix a $j$ such that $C$ meets each element $D\in\mathscr{E}_{j}$. Select a $D_{j}\in\mathscr{E}_{j}$. Then $D_{j}\subseteq\bigcup\mathscr{E}\subseteq E$ and hence $$\emptyset\neq C\cap D_{j}\subseteq C\cap E.$$ This leads to $C\in p_{E}$. So we have $\bigcup\limits_{j\in J}H_{j}\subseteq p_{E}$.

Suppose $G\in p_{E}$ or equivalently, $G\in\mathcal{O}(X)$ and $G\cap E\neq\emptyset$. This implies $G\cap(\bigcup\mathscr{E})\neq\emptyset$. Take an $E_{0}\in\mathscr{E}$ with $G\cap E_{0}\neq\emptyset$. In other words, $E_{0}\in\mathscr{E}\cap\diamondsuit G$ and $\mathscr{E}\setminus(\diamondsuit G)$ is a proper closed subset of $\mathscr{E}$. By the minimality of $\mathscr{E}$, there exists a $j_{1}$ such that $\mathscr{C}_{j_{1}}\cap(\mathscr{E}\setminus(\diamondsuit G))=\emptyset$. As a consequence, $\mathscr{E}_{j_{1}}=\mathscr{C}_{j_{1}}\cap\mathscr{E}\subseteq\diamondsuit G$ and therefore $G\in H_{j_{1}}$. Thus $p_{E}$ is contained in $\bigcup\limits_{j\in J}H_{j}$.

So we get $p_{E}=\bigcup\limits_{j\in J}H_{j}$. As $\mathscr{H}_{X}$ satisfies the second rule in Definition 5.4, $p_{E}\in\mathscr{H}_{X}$. Clearly, $E$ is nonempty and closed in $X$. This yields that $E\in\mathscr{P}$.By Lemma 5.3(3), $E$ is also irreducible in $X$. Note that $\mathscr{E}=cl_{\tau}(\mathscr{E})=\Box E=cl_{\tau}(\{E\})$. This means that $E\in\mathscr{E}\subseteq\mathscr{D}$. As $\mathscr{E}=cl_{\tau}(\{E\})$ meets each $\mathscr{C}_{j}$ and $\mathscr{C}_{j}$ is saturated, we deduce that $E$ belongs to each $\mathscr{C}_{j}$. As a result, $E\in\mathscr{D}\cap(\bigcap\limits_{j\in J})\mathscr{C}_{j}$. Therefore, $(\mathscr{P},\tau)$ is well-filtered.
\end{proof}

\begin{theorem}
{\rm Let $X$ be a $T_0$ space and $A\subseteq X$ a closed subset. Then
$A$ is a WD set of $X$ if and only if $p_A\in\mathscr H_X$.}
\end{theorem}

\begin{proof} Let $p_{A}\in\mathscr H_X$. Then by Lemma 5.3, $A$ is an irreducible closed subset of $X$. Suppose $Y$ is a well-filtered space and $f:X\longrightarrow Y$ is continuous. By Lemma 5.7 and the minimality of $\mathscr H_X$, $\mathscr H_X\subseteq R_{f}$. So $p_{A}\in R_{f}$ and hence there is a $K_{A}\in K(Y)$ such that for each $V\in\mathcal{O}(Y)$,
$$
\begin{aligned}
K_{A}\subseteq V
&\quad\Longleftrightarrow\quad f^{-1}(V)\in p_A\\
&\quad\Longleftrightarrow\quad f(A)\cap V\ne\varnothing.
\end{aligned}
$$
Suppose $\{V_j\mid j\in J\}\subseteq\mathcal{O}(Y)$ covers $K_{A}$. Then we have $f(A)\cap(\bigcup\limits_{i\in J}V_{j})\neq\emptyset$. So we can fix a $j_{0}\in J$ with $f(A)\cap V_{j_{0}}\neq\emptyset$. As a consequence, $K_{A}\subseteq V_{j_{0}}$. Hence $K_{A}$ is supercompact in $Y$. By Lemma 5.6, $K_{A}=\up_{Y} y$ for some $y\in Y$. Now, we have that for each $\widetilde{V}\in\mathcal{O}(Y)$, $$y\in \widetilde{V}\Longleftrightarrow K_{A}\subseteq\widetilde{V}\Longleftrightarrow f(A)\cap\widetilde{V}\neq\emptyset.$$ So we have $cl_{Y}f(A)=cl_{Y}\{y\}$. Therefore, $A$ is a WD set.

Conversely, let $A$ be a WD set of $X$. By Lemma 5.9, $(\mathscr{P},\tau)$ is well-filtered.
Define a mapping
$\eta:X\longrightarrow \mathscr{P}$ by $\eta(x)=cl_X(\{x\})$. Then for each $U\in\mathcal{O}(X)$, we have
$\eta^{-1}(\diamondsuit U)=\{x\in X\mid cl_X(\{x\})\cap U\neq\emptyset\}=U$. So $\eta$ is continuous. As $A$ is a WD set, there exists $C_0\in \mathscr{P}$ such that
$$
  cl_\tau(\eta(A))= cl_\tau(\{C_0\}).
$$
Then we obtain $$A=cl_{X}(A)=cl_{X}(\bigcup\eta(A))=\bigcup cl_\tau(\eta(A))=\bigcup cl_\tau(\{C_0\})=\bigcup\Box C_{0}=C_{0}.$$
As $C_{0}\in \mathscr{P}$, we know $p_A=p_{C_0}\in\mathscr H_X$.
The proof is complete.\end{proof}

\section{Further Work}

In the following, we list some problems that are no answers.

\begin{problem} {\rm Let $P$ and $Q$ be dcpos such that

 (1)$\Sigma K(P)$ and $\Sigma K(Q)$ are well-filtered;

 (2)$K(P\times Q)$ is a dcpo.

Is $\Sigma K(P\times Q)$ well-filtered? In particular, the other conditions keep unchanged and $Q$ is a continuous dcpo. Is $\Sigma K(P\times Q)$ well-filtered?}
\end{problem}

\begin{problem}{\rm Let $P$ and $Q$ be dcpos such that

 (1)$\Sigma K(P)$ and $\Sigma K(Q)$ are sober;

 (2)$K(P\times Q)$ is a dcpo.

 Is $\Sigma K(P\times Q)$ sober?}
\end{problem}


\section*{References}
\begin{thebibliography}{s9}

\bibitem{AJ94}
S. Abramsky, A. Jung, Domain theory, in: S. Abramsky, D. M. Gabbay
and T. S. E. Maibaum (Eds.), Handbook of Logic in Computer Science,
vol. 3, Clarendon Press, Oxford, 1994, pp. 1--168.

\bibitem{GG03} G. Gierz, K. H. Hofmann, K. Keimel, J. D. Lawson, M. W. Mislove and D. S. Scott, Continuous Lattices and Domains, Encyclopedia of Mathematics and its Applications, vol. 93, Cambridge University Press, Cambridge, 2003.

\bibitem{JG13} J. Goubault-Larrecq, Non-Hausdorff Topology and Domain Theory, New Mathematical Monographs, vol. 22, Cambridge University Press, New York, 2013.

\bibitem{Heckmann92}
R. Heckmann, An upper power domain construction in terms of strongly
compact sets, in: S. Brookes, M. Main, A. Melton, M. Mislove and
D. Schmidt (Eds.), Mathematical Foundations of Programming Semantics,
Lecture Notes in Computer Science, vol. 598, Springer, Berlin,
1992, pp. 272--293.

\bibitem{HK13}
R. Heckmann, K. Keimel, Quasicontinuous domains and the Smyth
powerdomain, Electronic Notes in Theoretical Computer Science,
298(2013)215--232.





\bibitem{HM81}
K. H. Hofmann, M. W. Mislove, Local compactness and continuous
lattices, in: B. Banaschewski and R.-E. Hoffmann (Eds.), Continuous
Lattices, Lecture Notes in Mathematics, vol. 871, Springer, Berlin,
1981, pp. 209--248.

\bibitem{MU13} R. Heckmann, K. Keimel, Quasicontinuous domains and the Smyth powerdomain, Electronic Notes in Theoretical Computer Science, 298 (2013) 215-232.

\bibitem{HMOOP281}Z. He, The product of two Scott sober complete lattices is not always Scott sober,arXiv:2609.39246, 2026.




\bibitem{Isbell82}
J. R. Isbell, Completion of a construction of Johnstone, Proceedings
of the American Mathematical Society, 85(1982)333--334.

\bibitem{Johnstone81}
P. T. Johnstone, Scott is not always sober, in: B. Banaschewski and
R.-E. Hoffmann (Eds.), Continuous Lattices, Lecture Notes in
Mathematics, vol. 871, Springer, Berlin, 1981, pp. 282--283.

\bibitem{E2} X. Jia, Meet-Continuity and Locally Compact Sober Dcpos, Ph.D Thesis, University of Birmingham, 2018.

\bibitem{Kou01}
H. Kou, $U_K$-admitting dcpo's need not be sober, in: Domains and
Processes, Kluwer Academic Publishers, Dordrecht, 2001, pp. 41--50.

\bibitem{LWX20}
J. D. Lawson, G. Wu and X. Xi, Well-filtered spaces, compactness,
and the lower topology, Houston Journal of Mathematics,
46(1)(2020)283--294.

\bibitem{MXJLZ26}
H. Miao, X. Xi, X. Jia, Q. Li and D. Zhao, Products of two sober
directed complete posets need not be sober, Journal of the London
Mathematical Society, 113(6)(2026)e70606.

\bibitem{Plotkin76}
G. D. Plotkin, A powerdomain construction, SIAM Journal on Computing,
5(3)(1976)452--487.

\bibitem{Rudin81}
M. E. Rudin, Directed sets which converge, in: L. F. McAuley and
M. M. Rao (Eds.), General Topology and Modern Analysis, Academic
Press, New York, 1981, pp. 305--307.

\bibitem{Scott72}
D. S. Scott, Continuous lattices, in: F. W. Lawvere (Ed.), Toposes,
Algebraic Geometry and Logic, Lecture Notes in Mathematics,
vol. 274, Springer, Berlin, 1972, pp. 97--136.

\bibitem{SXXZ19}
C. Shen, X. Xi, X. Xu and D. Zhao, On well-filtered reflections of
$T_0$ spaces, Topology and its Applications, 267(2019)106869.

\bibitem{Smyth78}
M. B. Smyth, Power domains, Journal of Computer and System Sciences,
16(1)(1978)23--36.

\bibitem{XL17}
X. Xi, J. D. Lawson, On well-filtered spaces and ordered sets,
Topology and its Applications, 228(2017)139--144.

\bibitem{XiZhao17}
X. Xi, D. Zhao, Well-filtered spaces and their dcpo models,
Mathematical Structures in Computer Science, 27(4)(2017)507--515.

\bibitem{XSXZ20}
X. Xu, C. Shen, X. Xi and D. Zhao, On $T_0$ spaces determined by
well-filtered spaces, Topology and its Applications,
282(2020)107323.

\bibitem{XWX23}
X. Xu, X. Wen and X. Xi, Scott topology on Smyth power posets,
Mathematical Structures in Computer Science, 33(9)(2023)832--867.

\bibitem{XZ21}
X. Xu, D. Zhao, Some open problems on well-filtered spaces and sober
spaces, Topology and its Applications, 301(2021)107540.









\end{thebibliography}
\end{document}